\input amstex
\magnification=1200
\documentstyle{amsppt}
\NoRunningHeads
\NoBlackBoxes
\define\Diff{\operatorname{Diff}}
\define\Rot{\operatorname{Rot}}
\define\PSL{\operatorname{PSL}}

\define\sltwo{\operatorname{\frak s\frak l}(2,\Bbb C)}
\define\ad{\operatorname{ad}}
\define\Hilb{\operatorname{Hilb}}
\define\form{\operatorname{form}}
\define\Pontr{\operatorname{Pontr}}
\define\vir{\operatorname{vir}}
\define\Vect{\operatorname{Vect}}
\define\sgn{\operatorname{sgn}}
\define\End{\operatorname{End}}
\define\Vir{\operatorname{Vir}}

\define\spn{\operatorname{span}}
\define\Gr{\operatorname{Gr}}
\define\Fg{\frak g}
\define\Fv{\frak v}
\define\Fw{\frak w}
\define\Fp{\frak p}
\define\CA{\Cal A}
\define\CF{\Cal F}
\define\CS{\Cal S}
\define\Der{\operatorname{Der}}
\topmatter
\title
Infinite dimensional geometry of $M_1\!=\!\Diff_+(\Bbb S^1)/\PSL(2,\Bbb
R)$ and $q_R$--conformal symmetries. II. Geometric quantization 
and hidden symmetries of Verma modules over Virasoro algebra.
\endtitle
\author Denis V. Juriev
\endauthor
\date July, 25, 1988\qquad\qquad\qquad\qquad math.RT/9807142
\enddate
\endtopmatter
\document
This article being a continuation of the first part [1] is addressed to 
specialists in representation theory, infinite dimensional geometry, 
quantum algebra, mathematical physics and informatics of interactive systems. 
The interrelations of the infinite dimensional geometry of the homogeneous 
K\"ahler manifold $M=\Diff_+(\Bbb S^1)/\Bbb S^1$ (see e.g. [2] and numerous 
refs wherein), its quotient $M_1=\Diff_+(\Bbb S^1)/\PSL(2,\Bbb R)$ and the 
$q_R$--conformal symmetries are discussed. Considering $q_R$--conformal 
symmetries in the Verma modules as Nomizu operators [1] one receives a 
connection in a bundle over $M_1$ with fibers isomorphic to Verma modules over 
the Lie algebra $\sltwo$. An interpretation of this connection in terms of 
geometric quantization is proposed. Simultaneously, another interpretation of 
$q_R$--conformal symmetries based on a geometric construction of the Verma 
modules over Virasoro algebra by the orbit method [3] as hidden symmetries 
is formulated. Such hidden symmetries are local over $M_1$, i.e. commute with 
the natural action of the commutative ring $\Cal O(M_1)$ of all polynomial
germs of holomorphic functions on $M_1$ in the Verma modules over Virasoro 
algebra. 

\head 1. Preliminary definitions\endhead

\subhead 1.1. The Lie algebra $\sltwo$ and Verma modules over it.
Lobachevski{\v\i}-Berezin $C^*$--algebra and $q_R$--conformal hidden
symmetries in the Verma modules over $\sltwo$\endsubhead
Lie algebra $\sltwo$ is a three-dimensional space of $2\times2$ complex
matrices with zero trace supplied by the standard commutator $[X,Y]=XY-YX$,
where the right hand side multiplication is the standard matrix multiplication.
In the basis
$$ l_{-1}=\left(\matrix0&1\\0&0\endmatrix\right),\quad
l_0=\left(\matrix\tfrac12&0\\0&-\tfrac12\endmatrix\right),\quad
l_1=\left(\matrix0&0\\-1&0\endmatrix\right) $$
the commutation relations are of the form: 
$$[l_i,l_j]=(i-j)l_{i+j}\qquad\qquad(i,j=-1,0,1).$$
Lie algebra $\sltwo$ is $\Bbb Z$--graded:
$\deg(l_i)=-\ad(l_0)l_i=i$, where $\ad(X)$ is the adjoint action operator
in the Lie algebra: $\ad(X)Y=[X,Y]$. Therefore, $\Bbb Z$--graded modules
over $\sltwo$ are $l_0$--diagonal. A vector $v$ in a $\Bbb Z$--graded
module over the Lie algebra $\sltwo$ is called extremal iff $l_1v=0$
and the linear span of vectors $l_{-1}^nv$ ($n\in\Bbb Z_+$) coincides
with the module itself (i.e. $v$ is a cyclic vector). A $\Bbb Z$--graded
module with an extremal vector (in this case it is defined up to a
multiplier) is called extremal [4]. An extremal module is called the
Verma module [5] iff the action of $l_{-1}$ is free in it, i.e. the vectors
$l_{-1}^nv$ are linearly independent. In the case of the Lie algebra $\sltwo$
the Verma modules are just the same as infinite dimensional extremal
modules. An extremal weight of the Verma module is the number defined
by the equality $l_0v=hv$, where $v$ is the extremal vector. The Verma
modules are defined for all complex numbers $h$ and are pairwise nonisomorphic.
Below we shall consider the Verma modules with real extremal weights only.

The Verma module $V_h$ over the Lie algebra $\sltwo$ with the extremal weight
$h$ may be realized in the space $\Bbb C[z]$ of polynomials of a complex
variable $z$. The formulas for the generators of the Lie algebra $\sltwo$ are
of the form:
$$l_{-1}=z,\quad l_0=z\partial_z+h,\quad l_1=z\partial_z^2+2h\partial_z,$$
here $\partial_z=\frac{d}{dz}$.

The Verma module is nondegenerate (i.e. does not contain any proper
sub\-mo\-dule) iff $h\ne-\tfrac{n}2$ ($n\in\Bbb Z_+$). The Verma module $V_h$
is called unitarizable (or Hermitian) iff it admits a structure of the
pre-Hilbert space such that $l_i^*=l_{-i}$. The completion of the
unitarizable Verma module will be denoted by $V^{\Hilb}_h$. The Lie algebra
$\sltwo$ acts in $V^{\Hilb}_h$ by the unbounded operators. Also it is rather
useful to consider the formal Verma modules $V^{\form}_h$, which are realized
in the space $\Bbb C[[z]]$ of formal power series of a complex variable $z$,
whereas the formulas for generators of the Lie algebra $\sltwo$ coincide
with ones above. Note that $V_h\subseteq V^{\Hilb}_h\subseteq V^{\form}_h$
and modules $V_h$, $V^{\Hilb}_h$, $V^{\form}_h$ form the Gelfand triple or
the Dirac equipment of the Hilbert space $V^{\Hilb}_h$. An action of the real
form of the Lie algebra $\sltwo$ generated by the anti-Hermitean operators
$il_0$, $l_1-l_{-1}$ and $i(l_1+l_{-1})$ in the
Hilbert space $V^{\Hilb}_h$ by the unbounded operators is exponentiated to
a unitary representation of the corresponding simply connected Lie group.

In the nonunitarizable Verma module over the Lie algebra $\sltwo$ there
exists the unique (up to a scalar multiple) indefinite sesquilinear
form $(\cdot,\cdot)$ such that $(l_iv_1,v_2)=(v_1,l_{-i}v_2)$ for any two
vectors $v_1$ and $v_2$ from the Verma module. If this sesquilinear form
is nondegenerate (in this case the Verma module is nondegenerate) then it
has a signature $(n,\infty)$, where $n$ is finite, and therefore, there is
defined a Pontryagin completion [6] of the Verma module. The corresponding
module in which the Lie algebra $\sltwo$ acts by the unbounded operators
will be denoted by $V^{\Pontr}_h$. The following chain of inclusions holds:
$V_h\subseteq V^{\Pontr}_h\subseteq V^{\form}_h$. An action of the real
form of the Lie algebra $\sltwo$ generated by the anti-Hermitean (with
respect to the nondegenerate indefinite sesquilinear form $(\cdot,\cdot)$)
operators $il_0$, $l_1-l_{-1}$ and $i(l_1+l_{-1})$ in the
Pontryagin space $V^{\Pontr}_h$ by the unbounded operators is exponentiated
to a pseudounitary representation of the corresponding simply connected
Lie group.

Let us describe some hidden symmetries in the Verma modules over the Lie
algebra $\sltwo$.

\proclaim{Proposition 1 [7]}\it In the nondegenerate Verma module $V_h$
over the Lie algebra $\sltwo$ there are uniquely defined the operators
$D$ and $F$ such that
$$\aligned [D,l_{-1}]=1,\quad [D,l_0]=D,\quad [D,l_1]=D^2,\\
[l_{-1},F]=1,\quad [l_0,F]=F,\quad [l_1,F]=F^2.\endaligned$$
If the Verma modules are realized in the space $\Bbb C[z]$ of polynomials
of a complex variable $z$ then
$$D=\partial_z,\qquad F=z\tfrac1{\xi+2h},$$
where $\xi=z\partial_z$. The operators $F$ and $D$ obey the following
relations:
$$[FD,DF]=0,\qquad [D,F]=q_R(1-DF)(1-FD),$$
where $q_R=\tfrac1{2h-1}$. In the unitarizable Verma module ($q_R\!\ne\!0$)
the operators $F$ and $D$ are bounded and $F^*=D$, $D^*=F$.
\endproclaim

\rm
The algebra generated by the variables $t$ and $t^*$ with the relations
$[tt^*,t^*t]=0$ and $[t,t^*]=q_R(1-tt^*)(1-t^*t)$ being the Berezin
quantization of the Lobachevski{\v\i} plane realized in the unit complex disc
(the Poincar\`e realization) [8] is called the Lobachevski{\v\i}-Berezin 
algebra. Proposition 1 allows to consider the Lobachevski{\v\i}-Berezin 
algebra as a $C^*$--algebra. 

\proclaim{Proposition 2 [7]}\it In the nongenerate Verma module $V_h$
over the Lie algebra $\sltwo$ there are uniquely defined the operators $L_n$
($n\in\Bbb Z$) such that
$$[l_i,L_n]=(i-n)L_{i+n}\qquad (i=-1,0,1;\ n\in\Bbb Z),$$
moreover, $L_i=l_i$ ($i=-1,0,1$). If the Verma modules are realized in the
space $\Bbb C[z]$ of polynomials of a complex variable $z$ then
$$L_k=(xi+(k+1)h)\partial_z^k\ (k\ge0),\quad
L_{-k}=z^k\tfrac{\xi+(k+1)h}{(\xi+2h)\ldots(\xi+2h+k-1)}\ (k\ge1),$$
where $\xi=z\partial_z$. The operators $L_n$ obey the following relations:
$$[L_n,L_m]=(n-m)L_{n+m},\text{\rm\ if}\ n,m\ge-1\text{\rm\ or}\
n,m\le1.$$
In the unitarizable Verma module the operators $L_n$ are unbounded and
$L_i^*=L_{-1}$.
\endproclaim

\rm
The operators $L_n$ are called the $q_R$--conformal symmetries. They may
be {\it sym\-bo\-li\-cal\-ly\/} represented in the form:
$$L_n=D^{nh}L_0D^{n(1-h)},\qquad L_{-n}=F^{n(1-h)}L_0F^{nh}.$$
To supply the symbolical recording by a sense one should use the general
com\-mu\-tation relations
$$[L_n,f(D)]=(-D)^{n+1}f'(D)\ (n\ge-1),\qquad
[L_{-n},f(F)]=F^{n+1}f'(D)\ (n\ge-1)$$
for $n=0$.

The commutation relations for the operators $D$, $F$ and the generators of
$q_R$--conformal symmetries with the generators of the Lie algebra $\sltwo$
mean that the families $J_k$ and $L_k$ ($k\in\Bbb Z$), where $J_i=D^i$,
$J_{-i}=F^i$ ($i\in\Bbb Z_+$), are families of tensor operators [6,9] for
the Lie algebra $\sltwo$.

\subhead 1.2. Infinite dimensional $\Bbb Z$--graded Lie algebras:
the Witt algebra $\Fw^{\Bbb C}$ of all Laurent polynomial vector
fields on a circle and the Virasoro algebra $\vir^{\Bbb C}$, its
one-dimensional nontrivial central extension. Group $\Diff_+(\Bbb S^1)$ 
of diffeomorphisms of a circle $\Bbb S^1$ and the Virasoro-Bott group $\Vir$. 
Flag manifold $M=\Diff_+(\Bbb S^1)/\Bbb S^1$ of the Virasoro-Bott group and 
infinite dimensional K\"ahler manifold $M_1=\Diff_+(\Bbb S^1)/\PSL(2,\Bbb R)$
\endsubhead
The Lie algebra $\Vect(S^1)$ is realized in the space of $C^\infty$--smooth
vector fields $v(t)\partial_t$ on a circle $S^1\simeq\Bbb R/2\pi\Bbb Z$
with the commutator
$$[v_1(t)\partial_t,v_2(t)\partial_t]=
(v_1(t)v_2'(t)-v_1'(t)v_2(t))\partial_t.$$
In the basis 
$$s_n=\sin(nt)\partial_t,\qquad c_n=\cos(nt)\partial_t,\qquad h=\partial_t$$
the commutation relations have the form:
$$\aligned
[s_n,s_m]&=\tfrac12((m-n)s_{n+m}+\sgn(n-m)(n+m)s_{|n-m|}),\\
[c_n,c_m]&=\tfrac12((n-m)s_{n+m}+\sgn(n-m)(n+m)s_{|n-m|}),\\
[s_n,c_m]&=\tfrac12((m-n)c_{n+m}-(m+n)c_{|n-m|})-n\delta_{nm}h,\\
[h,s_n]&=nc_n,\quad [h,c_n]=-ns_n.
\endaligned$$
Let us denote by $\Vect^{\Bbb C}(S^1)$ the complexification of the Lie algebra
$\Vect(S^1)$. In the basis $e_n=ie^{ikt}\partial_t$ the commutation relations
in the Lie algebra Li $\Vect^{\Bbb C}(S^1)$ have the form:
$$[e_j,e_k]=(j-k)e_{j+k}.$$

It is rather convenient to consider an imbedding of the circle $S^1$ into
the complex plane $\Bbb C$ with the coordinate $z$, so that $z=e^{it}$ on the
circle and the elements of the basis $e_k$ ($k\in\Bbb Z$) are represented by
the Laurent polynomial vector fields: $$e_k=z^{k+1}\partial_z.$$ 
The $\Bbb Z$--graded Lie algebra generated by the Laurent polynomial vector
fields (i.e. by the finite linear combinations of elements of the basis
$e_k$) is called the Witt algebra and is denoted by $\Fw^{\Bbb C}$.
The Witt algebra $\Fw^{\Bbb C}$ is the complexification of the
subalgebra $\Fw$ of the algebra $\Vect(S^1)$ generated by the
trigonometric polynomial vector fields on a circle $S^1$, i.e. by the finite
linear combinations of elements of the basis $s_n$, $c_n$ and $h$.

The Lie algebra $\Vect(S^1)$ admits a nontrivial one-dimensional central
extension defined by the Gelfand-Fuchs 2-cocycle [10]:
$$c(v_1(t)\partial_t,v_2(t)\partial_t)=
\int_0^{2\pi}(v_1'(t)v_2''(t)-v_2'(t)v_1''(t))\,dt.$$
This extension being continued to the complexification $\Vect^{\Bbb C}(S^1)$
of the Lie algebra $\Vect(S^1)$ and reduced to the subalgebra
$\Fw^{\Bbb C}$ defines a central one-dimensional ex\-tension of the Witt
algebra, which is called the Virasoro algebra and is denoted by
$\vir^{\Bbb C}$. The Virasoro algebra is generated by the elements
$e_k$ ($k\in\Bbb Z$) and the central element $\bold c$ with the commutation
relations:
$$[e_j,e_k]=(j-k)e_{j+k}+\frac{j^3-j}{12}\bold c\delta(j+k)$$
and is the complexification of a central extension $\vir$ of the Lie algebra
$\Fw$. 
Because the 2-cocycle
$$\int_0^{2\pi}(v_1(t)v_2'(t)-v_2(t)v_1'(t))\,dt$$
is trivial, the cocycle, which defines the Virasoro algebra, is indeed 
equivalent to the Gelfand-Fuchs 2-cocycle above and is known under the same 
name. The modified Gelfand-Fuchs cocycle is $\sltwo$--invariant, where 
$\sltwo$ is a subalgebra of $\vir^{\Bbb C}$ generated by $e_{-1}$, $e_0$ and
$e_1$, so it is handier in practice. Below we shall use the mo\-di\-fied 
version of Gelfand-Fuchs 2-cocycle under this name only. In the irreducible 
representation the central element $\bold c$ of the Virasoro algebra is mapped 
to a scalar operator, which is proportional to the identity operator with a
coefficient $c$ called {\it the central charge}.

Let $\Diff(\Bbb S^1)$ denote the group of all diffeomorphisms of the circle
$\Bbb S^1)$. The group manifold $\Diff(\Bbb S^1)$ splits into two 
connected components, the subgroup $\Diff_+(\Bbb S^1)$ and the coset
$\Diff_-(\Bbb S^1)$. The diffeomorphisms in $\Diff_+(\Bbb S^1)$ preserve
the orientation on the circle $\Bbb S^1)$ and those in $\Diff_-(\Bbb S^1)$
reverse it. The Lie algebra of $\Diff_+(\Bbb S^1)$ can be identified with
$\Vect(\Bbb S^1)$.

The infinite-dimensional group $\Vir$ corresponding to the Virasoro algebra
$\vir$ (more precisely, to the central extension of the Lie algebra
$\Vect(\Bbb S^1)$ defined by the Gelfand-Fuchs cocycle, whereas the
Virasoro algebra $\vir$ is an extension of the real form $\Fw$ of the
Witt algebra $\Fw^{\Bbb C}$) is a central extension of the group
$\Diff(\Bbb S^1)$. The corresponding 2-cocycle was calculated by R.Bott [11]:
$$c(g_1,g_2)=\int\log(g_1'\circ g_2)\,\log(g_2').$$
The group $\Vir$ is called the Virasoro-Bott group.

The flag manifold $M$ of the Virasoro-Bott group is a homogeneous space
with transformation group $\Diff_+(\Bbb S^1)$ and isotropy group $\Bbb S^1)$.
There exist several different realizations of this manifold [12-15].

{\it Algebraic realization}. The space $M$ can be realized as a conjugacy 
class in the group $\Diff_+(\Bbb S^1)$ or in the Virasoro-Bott group $\Vir$.

{\it Probabilistic realization}. Let $P$ be the space of real probability
measures $\mu=u(t)\,dt$ with smooth positive density $u(t)$ on $\Bbb S^1$.
The group $\Diff_+(\Bbb S^1)$ naturally acts on $P$ by the formula
$$g:u(t)\,dt\mapsto u(g^{-1}(t))\,dg^{-1}(t).$$
The action is transitive and the stabilizer of the point $(2\pi)^{-1}\,dt$
is isomorphic to $\Bbb S^1$, therefore, $P$ can be identified with $M$.

{\it Orbital realization}. The space $M$ can be considered as an orbit of the
coadjoint representation of $\Diff_+(\Bbb S^1)$ or $\Vir$. Namely, the
elements of the dual space $\vir^*$ of the Virasoro algebra $\vir$ are
identified with the pairs $(p(t)\,dt^2,b)$; the coadjoint action of $\Vir$
has the form
$$K(g)(p,b)=(gp-bS(g),b),$$
where
$$S(g)=\frac{g'''}{g'}-\frac32\left(\frac{g''}{g'}\right)^2$$
is the Schwarzian (the Schwarz derivative) and $gp$ denotes the natural action
of $g$ on the quadratic differential $p$. The orbit of the point $(a\cdot
dt^2,b)$ coincides with $M$ if $a/b\ne-n^2/2$, $n=1,2,3,\cdots$. Therefore,
a family $\omega_{a,b}$ of symplectic structures (Kirillov forms) is defined
on $M$ (cf.[16]).

{\it Analytic realization}. Let us consider the space $S$ of univalent
functions $f(z)$ on the unit disk $D_+$ such that $f(0)=0$, $f'(0)=1$ and
$f'(e^{it})\ne0$ [17-19] (ordinarily the least condition is omitted;
so below we shall consider a proper subclass of a conventional class
$S$ under the same notation). The Taylor coefficients $c_1, c_2, c_3,\cdots$ 
in the expansion
$$f(z)=z+c_1z^2+c_2z^3+\cdots+c_nz^{n+1}+\cdots$$
form a coordinate system on $S$. The class $S$ can be naturally identified
with $M$ via the Kirillov construction [14]. The Lie algebra $\Vect(\Bbb S^1)$
acts on $S$ by the formulas 
$$L_vf(z)=-if^2(z)\oint\left(\frac{wf'(w)}{f(w)}\right)^2\frac{v(w)}
{f(w)-f(z)}\frac{dw}w.$$
The Kirillov construction supply $M$ by the complex structure. The symplectic
structure $\omega_{a,b}$ coupled with the complex structure determines a
K\"ahler metric $w_{a,b}$ on $M$. More detailed information on the 
infinite-dimensional geometry of the flag manifold $M$ is contained in [2]
(see also refs wherein). Note only that the curvature tensors of the K\"ahler
connections on $M$ were calculated in [15].

The subgroup $\Bbb S^1$ is contained in each of the subgroups $H_k$,
$k=1,2,3,\cdots$, generated by the generators $ie_0$, $e_k-e_{-k}$ and
$i(e_k+e_{-k})$. The subgroup $H_k$ is isomorphic to the $k$-folded
covering of the group $H_1=\PSL(2,\Bbb R)$ and acts on $\Bbb S^1$ by the
formulas
$$z\mapsto\left(\frac{az^k+b}{\bar bz^k+\bar a}\right)^{1/k}=
\alpha z\left(\frac{1+bz^{-k}}{1+\bar bz^k}\right)^{1/k},$$
where $\alpha=(a/|a|)^{2/k}$ is an univalued function on $H_k$.

The homogeneous space $M_k=\Diff_+(\Bbb S^1)/H_k$ is a symplectic manifold
and can be identified with the orbit with $a/b=-k^2/2$ in the coadjoint
representation of the Virasoro-Bott group. The manifold $M_k$ is a quotient
of $M$ by the $\Diff_+(\Bbb S^1)$-invariant foliation $\Cal F_k$ generated 
by $s_k$ and $c_k$ (considered as elements of a reductive basis on $M$ [15]).
Note that though $M$ is a reductive space the manifolds $M_k$ are not 
reductive. The foliations $\Cal F_k$ are complex foliations, which were
discovered by V.Yu.Ovsienko and O.D.Ovsienko. The foliation $\Cal F_1$
is holomorphic whereas other foliations are not holomorphic. Hence,
the almost complex structure on $M_1$ is integrable and, therefore,
$M_1=\Diff_+(\Bbb S^1)/\PSL(2,\Bbb R)$ is a homogeneous K\"ahler manifold.
This is not true for $k\ge2$. The main characteristics of the K\"ahler metrics 
on $M_1$ were calculated in [15] (really, there were calculated their 
back-liftings to $M$). In particular, these K\"ahler metrics are einsteinian. 

Let us consider any point $f(z)$ of $M$. A fiber of the first Ovsienko
foliation $\Cal F_1$ consists of points 
$$f_b(z)=\frac{f(z)}{1-bf(z)},$$
where $b\in\bar\Bbb C\setminus[f(D_+)]^{-1}$. Therefore, $M_1$ may be
identified with the class $S^{(1)}$ of the univalent functions $f(z)$ such 
that $f(0)=0$, $f'(0)=1$, $f''(0)=0$, $f'(e^{it})\ne0$. The variational
formulas for an action of $\Vect(\Bbb S^1)$ on $S^{(1)}$ are of the form
$$L^{(1)}_vf(z)=L_vf(z)-[L_vf]_2f^2(z),$$
where 
$$[L_vf]_2=\frac12(L_vf)''(0)=i\oint\frac{w(f'(w))^2v(w)\,dw}{f^3(w)}.$$

\subhead 1.3. Verma modules over the Virasoro algebra and their construction
by the orbit method [3]  
\endsubhead 
The algebras $\Fw^{\Bbb C}$ and $\vir^{\Bbb C}$ as a $\Bbb Z$-graded 
algebras possess Verma modules over them, which was studied by many authors
(see e.g. [20-22] and also [23,24]). Namely, put $\vir^{\Bbb C}_+=\spn(e_k,
k\ge0)$, $\chi_{h,c}$ be a character of $\vir^{\Bbb C}_+$ defined as 
$$\chi_{h,c}(e_k)=0,\text{\ if\ } k>0,\quad
\chi_{h,c}(e_0)=h,\quad \chi_{h,c}(c)=c.$$
Verma module $V_{h,c}$ is a $\vir^{\Bbb C}$--module induced from the
character $\chi_{h,c}$ of the subalgebra $\vir^{\Bbb C}_+$. Otherwords,
$$V_{h,c}=U(\vir^{\Bbb C})\otimes_{U(\vir^{\Bbb C}_+)}V_{\chi_{h,c}},$$
where $V_\chi$ is a $\vir^{\Bbb C}_+$--module defined by a character $\chi$.
$U(\vir^{\Bbb C})$ and $U(\vir^{\Bbb C}_+)$ are the universal envelopping
algebras of the Lie algebras $\vir^{\Bbb C}$ and $\vir^{\Bbb C}_+$.
The Verma module $V_{h,c}$ is $\Bbb Z$--graded $\vir^{\Bbb C}$--module;
if $h$ and $c$ are real (what will be supposed below) there is defined
the unique up to a multiple invariant Hermitian form in $V_{h,c}$. The
Verma module is unitarizable if and only if the Hermitian form is positive
definite; let us denote by $D_n(h,c)$ the determinant of this form in the
$n$-th homogeneous component of $V_{h,c}$ in the basis $e_1^{k_1}\ldots
e_j^{k_j}v$, $k_j\ge0$ ($v$ is the extremal vector, i.e. $e_0v=hv$,
$e_{-m}v=0$), then as it was shown by V.G.Kac, B.L.Feigin and D.B.Fuchs
$$D_n(h,c)=A\prod_{\gathered 0<\alpha\le\beta\\ \alpha\beta\le n\endgathered}
\Phi^{p(n-\alpha\beta)}_{\alpha,\beta},$$
where
$$\gathered
\Phi_{\alpha,\beta}(h,c)\!=\!(h\!+\!\frac{c\!-\!13}{24}(\beta^2\!-\!1)
\!+\!\frac{\alpha\beta\!-\!1}2)(h\!+\!\frac{c\!-\!13}{24}(\alpha^2\!-\!1)
\!+\!\frac{\alpha\beta\!-\!1}2)+\frac{\alpha^2\!-\!\beta^2}{16},\\
\Phi_{\alpha,\alpha}(h,c)=h+\frac{c-13}{24}(\alpha^2-1).
\endgathered
$$
If for any $\alpha$, $\beta$ $\Phi_{\alpha,\beta}(h,c)\ne0$, then the module 
$V_{h,c}$ is irreducible and is not contained in any other Verma module; if
there exists exactly one pair $(\alpha,\beta)$ such that $\Phi_{\alpha,\beta}
(h,c)=0$ then there three possibilities may be realized: 1) $\alpha\beta<0$,
then $V_{h,c}$ may be imbedded into the Verma module $V_{h+\alpha\beta,c}$,
2) $\alpha\beta>0$, then $V_{h,c}$ contains a submodule $V_{h+\alpha\beta,c}$,
3) either $\alpha=0$ or $\beta=0$, then $V_{h,c}$ is irreducible and is not
a submodule of another Verma module. If there exists two pairs $(\alpha_1,
\beta_1)$ and $(\alpha_2,\beta_2)$ such that $\Phi_{\alpha_i,\beta_i}(h,c)=0$,
then there exists an infinite number of pairs $(\alpha,\beta)$, which possess
such property; this situation is realized if
$$
\gathered
c_{1,2}=1-\frac{6((\alpha_1\pm\alpha_2)-(\beta_1\pm\beta_2))^2}
{(\alpha_1\pm\alpha_2)(\beta_1\pm\beta_2)},\\
h_{1,2}=\frac{(\alpha_2\beta_1-\alpha_1\beta_2)^2-((\alpha_1\pm\alpha_2)-
(\beta_1\pm\beta_2))^2}{4(\alpha_1\pm\alpha_2)(\beta_1\pm\beta_2)}.
\endgathered
$$
In this case the structure of Verma modules is described by the Feigin--Fuchs
theory.

The Verma module $V_{h,c}$ is unitarizable if $h>0$, $c>1$; the Verma module
$V_{h,c}$ contains an unitarizable quotient if (a) $h>0$, $c>1$;
(b) $c=1-\frac6{p(p+1)}$, $h=\frac{(\alpha p-\beta(p+1))^2-1}{4p(p+1)}$,
$\alpha,\beta,p\in\Bbb Z$; $p\ge2$; $1\le\alpha\le p$, $1\le\beta\le p-1$.

Let us now describe a geometric way of the construction of the Verma modules 
over the Virasoro algebra based on the orbit method [3]. It uses the following
facts:

(1) to any $\Diff_+(\Bbb S^1)$--invariant K\"ahler metric $w_{h,c}$ 
(new parameters $h$ and $c$ are related to $a$ and $b$ as $a=h-\frac{c}{12}$, 
$b=\frac{c}{12}$ and $h=a+b$, $c=12b$) on the space 
$M=\Diff_+(\Bbb S^1)/\Bbb S^1$ one can assign a linear holomorphic bundle 
$E_{h,c}$ over $M$ with the following properties:
\roster
\item"(a)" $E_{h,c}$ is the Hermitian bundle, with metric $\exp(-U_{h,c})
d\lambda d\bar\lambda$, where $\lambda$ is the coordinate in the fiber and
$K_{h,c}=\exp(U_{h,c})$ is the Bergman kernfunction, i.e. the exponential
of the K\"ahler potential of the metric $w_{h,c}$ (the K\"ahler potentials
for $w_{h,c}$ were calculated in [15,3]),
\item"(b)" the algebra $\vir^{\Bbb C}$ holomorphically acts in the prescribed
bundle by covariant derivatives with respect to the Hermitian connection with
curvature form $2\pi i\omega_{h,c}$ (here $\omega_{h,c}$ is the symplectic
structure corresponded to $w_{h,c}$);
\endroster

(2) let $\Cal O(E_{h,c})$ be the space of all polynomial germs (in some 
natural trivialization) of sections of the bundle $E_{h,c}$ in the coordinates
$c_1,\ldots c_k,\ldots$. The action of $\vir^{\Bbb C}$ in its $\Bbb Z$-graded
module $\Cal O(E_{h,c})$ ($\deg(c_k)=k$) is defined by the formulas
$$\aligned
L_p&=\partial_p+\sum_{k\ge1}(k+1)c_k\partial_{k+p}\quad (p>0),\qquad
L_0=\sum_{k\ge1}kc_k\partial_k+h,\\
L_{-1}&=\sum_{k\ge1}((k+2)c_{k+1}-2c_1c_k)\partial_k+2hc_1,\\
L_{-2}&=\sum_{k\ge 1}((k+3)c_{k+2}-(4c_2-c_1^2)c_k-b_k(c_1,\ldots c_{k+2}))
\partial_k\\
&\qquad\qquad+h(4c_2-c_1^2)+\frac{c}2(c_2-c_1^2),\\
L_{-n}&=\frac1{(n-2)!}\ad^{n-2}L_{-1}\cdot L_{-2}\quad (n>0),
\endaligned
$$ 
where $\partial_k=\frac{\partial}{\partial c_k}$ and $b_k$ are the Laurent
coefficients of the function $1/(zf(z))$.

Let us choose the basis $e^{a_1,\ldots,a_n}=c_1^{a_1}\ldots c_n^{a_n}$ in
$\Cal O(E_{h,c}^*)$ and let $\Cal O^*(E^*_{h,c})$ be the space of all linear
functionals $p$ on $\Cal O(E^*_{h,c})$ such that the relation $p(x)\ne0$
implies $\deg(x)\le N$. The space $\Cal O^*(E_{h,c}^*)$ is called the
Fock space of the pair $(M,E_{h,c})$ and is denoted by $F(E_{h,c})$.
The Verma module $V_{h,c}$ over $\vir^{\Bbb C}$ is realized in the Fock
space $F(E_{h,c})$. The formal Verma module $V^{\form}_{h,c}$ may be
obtained if one omits a condition on the functionals $p$. If we choose the
basis $e_{a_1,\ldots,a_n}=:c_1^{a_1}\ldots c_n^{a_n}$ in $F(E_{h,c}$ such that
$$\left<e_{a_1,\ldots,a_n},e^{b_1,\ldots,b_m}\right>=a_1!\ldots a_n!
\delta_n^m\delta_{a_1}^{b_1}\ldots\delta_{a_n}^{b_n},$$
then the action of the Virasoro algebra in this basis is given by the formulas
$$\aligned
L_{-p}&=c_p+\sum_{k\ge1}c_{k+p}\partial_k\quad(p>0),\qquad
L_0=\sum_{k\ge1}kc_k\partial_k+h,\\
L_1&=\sum_{k\ge1}c_k((k+2)\partial_{k+1}-2\partial_1\partial_k)+2h\partial_1,\\
L_2&=\sum_{k\ge1}c_k((k+3)\partial_{k+2}-(4\partial_2-\partial^2_1)\partial_k-
b_k(\partial_1,\ldots,\partial_{k+2}))\\
&\qquad\qquad+h(4\partial_2-\partial^2_1)+\frac{c}2(\partial_2-\partial_1^2),\\
L_n&=\frac{(-1)^n}{(n-2)!}\ad^{n-2}L_1\cdot L_2\quad (n>2).
\endaligned$$
The action of $\vir^{\Bbb C}$ is exponentiated to the projective action of
$\Vir$ in a certain subspace of $V^{\form}_{h,c}$. The Bergman kernfunction
expanded by $c_k$ and $\bar c_k$ coincides with the invariant Hermitian 
form in the Verma module $V_{h,c}$ written in the basis $e_{a_1,\ldots,a_n}$.
If the Verma module is unitarizable then the action of $\vir^{\Bbb C}$ is
exponentiated to the unitary action of the universal covering $\widetilde\Vir$
of the Virasoro-Bott group in the space $V^{\Hilb}_{h,c}$. In the 
nonunitarizable case the situation is slightly more complicated [23].

The construction above may be applied to the manifold $M_1$ instead of
$M$. Thr K\"ahler metrics $w_c$ on $M_1$ are just the reduction of 
degenerate metrics $w_{h,c}$ on $M$ with $h=0$. So the bundle $E_c$
over $M_1$ related to the metric $w_c$ is a restriction of $E_{0,c}$
over $M$ on $M_1$ (here $M_1$ is imbed into $M$ as a subclass $S^{(1)}$ of
$S$). The action of $\vir^{\Bbb C}$ in the space $\Cal O(E_c)$ is
defined by the formulas
$$\allowdisplaybreaks\align
L_p&=\partial_p+\sum_{k\ge0}(k+1)c_k\partial_{k+p}\quad(p\ge2),\\
L_1&=\sum_{k\ge2}(k+1)c_k\partial_{k+1}-\Gamma,\\
L_0&=\sum_{k\ge2}c_k\partial_k,\\
L_1&=\sum_{k\ge2}(k+2)c_{k+1}\partial_k-3c_2\Gamma,\\
L_2&=\sum_{k\ge2}((k+3)c_{k+2}-4c_2c_k-b_k(0,c_2,\ldots,c_{k+2}))
\partial_k+\frac{c}2c_2-5c_3\Gamma,\\
L_{-n}&=\frac1{(n-2)!}\ad^{n-2}L_{-1}\cdot L_{-2}\quad (n>0),
\endalign
$$
where
$$
\Gamma=2\sum_{k\ge2}c_k\partial_{k+1}+\sum_{i,j\ge2}c_ic_j\partial_{i+j+1}.
$$
The Fock space $F(E_c)$ is a $\Bbb Z$--graded $\vir^{\Bbb C}$--module
in which the module $W_c$ (a quotient of $V_{0,c}$) is realized.
The Virasoro algebra $\vir^{\Bbb C}$ acts in $W_c$ as
$$\align
L_{-p}&=c_p+\sum_{k\ge2}(k+1)c_{k+p}\partial_k\quad(p\ge2),\\
L_{-1}&=\sum_{k\ge2}(k+1)c_{k+1}\partial_k-\Gamma^*,\\
L_0&=\sum_{k\ge2}c_k\partial_k,\\
L_1&=\sum_{k\ge2}(k+2)c_k\partial_{k+1}-3\Gamma^*\partial_2,\\
L_2&=\sum_{k\ge2}c_k((k+3)\partial_{k+2}-4\partial_2\partial_k-
b_k(0,\partial_2,\ldots,\partial_{k+2})+\frac{c}2\partial_2-
5\Gamma^*\partial_3,\\
L_n&=\frac{(-1)^n}{(n-2)!}\ad^{n-2}L_1\cdot L_2\quad (n>2).
\endalign
$$
where
$$\Gamma^*=2\sum_{k\ge2}c_{k+1}\partial_k+\sum_{i,j\ge2}c_{i+j+1}\partial_i
\partial_j.$$

\head 2. Geometric quantization and hidden symmetries of Verma modules
over the Virasoro algebra\endhead

\subhead 2.1. Verma modules over the Virasoro algebra and geometric
quantization\endsubhead
Let us consider an arbitrary (topologically trivial for simplicity) K\"ahler 
manifold $N$ and a linear bundle $E$ over it. There exists a canonical 
imbedding of the manifold $N$ into the projectivization $\Bbb P(F(E))$ of 
the Fock space $F(E)$ of the bundle $E$. If $E$ is the Hermitian bundle and 
the scalar products in fibers is defined via the Bergman kernfunction (the 
exponent of a K\"ahler potential) of the K\"ahler metric $w$ on $N$ then $F(E)$ 
is a pseudo-Hilbert space and the imbedding $N\mapsto\Bbb P(F(E))$ is 
isometric. The bundle $E$ over $N$ may be restored as a restriction 
$\left.\tau\right|_N$ of the canonical bundle $\tau$ over $\Bbb P(F(E))$
onto the submanifold $N$. Note that $\tau$ possess a natural Hermitian
connection, which restriction onto $E=\left.\tau\right|_N$ is just the
prequantization connection in $E$ (the Hermitian connection is called a 
prequantization one if its curvature is equal to $2\pi\omega$, where the 
symplectic form $\omega$ is the imaginary part of the K\"ahler metric $w$ 
[25,26]). 

If the K\"ahler manifold $N$ is not topologically trivial then one is able
to consider it locally and then to globalize all constructions. At least,
any isometric imbedding of a K\"ahler manifold $N$ into a projective space
$\Bbb P(V)$ supplies $N$ by a prequantization bundle $E=\left.\tau\right|_N$
(a Hermitian bundle with a prequantization connection) over it.

A natural generalization of projective spaces are Grassmannians. It is 
rather interesting to generalize the geometric quantization schemes exposed
above onto imbeddings into Grassmannians. Indeed, let us consider a
K\"ahler manifold $N$ and its isometric imbedding into the Grassmannian
$\Gr(V;V_0)$ of all subspaces of the Hilbert (Euclidean) space $V$, which 
are isometric to its subspace $V_0$. One is able to consider a canonical 
bundle $\tau$ over $\Gr(V;V_0)$ with fibers isometric to $V_0$ and a 
canonical Hermitian connection $\nabla$ in it and then to restrict $\tau$ 
and $\nabla$ onto the submanifold $N$. Thus, one obtains a Hermitian bundle 
$\hat E$ over $N$ and a Hermitian connection $\nabla$ in it. In general, 
$\nabla$ is not flat and has a non-trivial (and nonscalar) curvature.

In particular, one may consider the K\"ahler manifold $N$ with a 
prequantization linear bundle $E$ and put $V=F(E)$. Let us also consider
a holomorphic bundle $N\mapsto N_1$. For any fiber $p$ of this bundle
let us put $V_0=F(p)\subseteq V$. Then, one receives an imbedding of
$N_1$ into $\Gr(V;V_0)$ and the construction above supplies $N_1$ by
the Hermitian bundle $\hat E$ over it with fibers isometric to $V_0$ and
Hermitian connection $\nabla$ in it.

One may specify the homogeneous case $N=G/H$, where $G$ is a real reductive
Lie group and $H$ is its certain subgroup, $E$ is a prequantization bundle 
over $N$. The Fock space $F(E)$ realizes a Verma module over $\Fg^{\Bbb 
C}$ [27]. Let also $G_0$ be a subgroup of $G$ such that $H\subseteq G_0$,
and $G_0/H\mapsto G/H\mapsto G/G_0$ is a holomorphic bundle. Then one 
obtains an imbedding $G/G_0\mapsto\Gr(V;V_0)$, where $V=F(E)$ and
$V_0=F(\left.E\right|_{G_0/H})$, the Hermitian bundle $\hat E$ over $G/G_0$
and Hermitian connection $\nabla$ in it. The situation is straightforwardly
generalized to infinite dimensional case (e.g. $G$ is the Virasoro-Bott group).

In the homogeneous case the bundle $\hat E$ possesses a natural action of
the Lie algebra $\Fg^{\Bbb C}$ by operators $L_X$ ($X\in\Fg^{\Bbb C}$). 
The covariant derivatives $\nabla_X$ along the elements of $\Fg^{\Bbb C}$
are not closed. In the physical terms their commutators realize the so-called
nonscalar Schwinger terms [28]. The operators $A_X=L_X-\nabla_X$ are fiberwise
and are called the Nomizu operators [29] of the connection $\nabla$.

The Fock space $F(\hat E)$ of the bundle $\hat E$ over $G/G_0$ is naturally
identified with the Fock space $F(E)$ of the bundle $E$ over $G/H$ and,
therefore, realizes a Verma module over $\Fg^{\Bbb C}$. Hence, the
covariant derivatives $\nabla_X$ and the Nomizu operators $A_X$ act in
the Verma module.

Specifying the situation for the Virasoro-Bott group one obtains the
following proposition.

\proclaim{Theorem 1} The Verma module $V_{h,c}$ over the Virasoro algebra
may be realized in the Fock space $F(\hat E_{h,c})$ of the vector bundle $\hat 
E_h$ over $M_1=\Diff_+(\Bbb S^1)/\PSL(2,\Bbb R)$ with fibers isomorphic to
the Verma module $V_h$ over the Lie algebra $\sltwo$. The Nomizu operators
$A_X$ ($X\in\vir^{\Bbb C}$) of the Hermitian connection $\nabla$ in $\hat E$ 
being considered in any fixed fiber coincide with the $q_R$--conformal 
symmetries in $V_h$.
\endproclaim

The theorem follows from the fact that the Nomizu operators $A_X$ are
tensor ones in $V_h$ with respect to $\sltwo$ and, therefore, coincide with 
$q_R$--conformal symmetries up to a multiple. However, for $X\in\sltwo$
$A_X$ are just the generators of $\sltwo$, hence, the multiple is equal to 1
\qed

This theorem is parallel to the main theorem of [1]. Below we shall discuss
the action of the operators $\nabla_X$ and $A_X$ ($X\in\vir^{\Bbb C}$) in the 
Verma modules over the Virasoro algebra.

\remark{Remark 1} $\hat E_{h,c}=\hat E_{h,0}\otimes E_c$, where $E_c$ is
a linear bundle over $M_1$. This representation is valid for the connections
$\nabla$ also.
\endremark

\subhead 2.2. Nomizu operators $A_X$ as hidden symmetries in the Verma 
modules over the Virasoro algebra\endsubhead
Let us realize the Verma modules $V_{h,c}$ over the Virasoro algebra
$\vir^{\Bbb C}$ in the Fock spaces $F(E_{h,c})$ of the Hermitian
bundles $E_{h,c}$. Then the commutative algebra $\Cal O(M_1)$ of all
polynomial germs of holomorphic functions on $M_1$ naturally acts in
$F(E_{h,c})$ as a subalgebra of $\Cal O(M)$.

\proclaim{Theorem 2} There exists a unique family of tensor operators
in the Verma module $V_{h,c}$ over the Virasoro algebra $\vir^{\Bbb C}$, 
which transform as the adjoint representation of $\vir^{\Bbb C}$, and
which are local over $M_1$, i.e. commute with the action of $\Cal O(M_1)$
in the module $V_{h,c}$. These tensor operators coincide with the
Nomizu operators $A_X$ up to a multiple.
\endproclaim 

Note that tensor operators from a family form the same family
being multiplied on a number.

A proof of the theorem 2 is just the same as of theorem 1.

\remark{Remark 2} 
$$[L_X,A_Y]=A_{[X,Y]}$$
for any $X,Y\in\vir^{\Bbb C}$.
\endremark

\subhead 2.3. Setting hidden symmetries free\endsubhead

\definition{Definition 1 [30]}

{\bf A.} A linear space $\Fv$ is called {\it a Lie composite\/}
iff there are fixed its subspaces $\Fv_1,\ldots \Fv_n$ ($\dim\Fv_i>1$)
supplied by the compatible structures of Lie algebras. Com\-pa\-ti\-bi\-li\-ty 
means that the structures of the Lie algebras induced in $\Fv_i\cap\Fv_j$ 
from $\Fv_i$ and $\Fv_j$ are the same. The Lie composite is called {\it 
dense\/} iff $\Fv_1\uplus\ldots\uplus\Fv_n=\Fv$ (here $\uplus$ denotes the sum 
of linear spaces). The Lie composite is called {\it connected\/}
iff for all $i$ and $j$ there exists a sequence $k_1,\ldots k_m$ ($k_1=i$,
$k_m=j$) such that $\Fv_{k_l}\cap\Fv_{k_{l+1}}\ne\varnothing$.

{\bf B.} {\it A representation\/} of the Lie composite $\Fv$ in the space $H$
is the linear mapping $T:\Fv\mapsto\End(H)$ such that $\left.T\right|_{\Fv_i}$
is a representation of the Lie algebra $\Fv_i$ for all $i$.

{\bf C.} Let $\Fg$ be a Lie algebra. A linear mapping $T:\Fg\mapsto\End(H)$
is called {\it the composed representation\/} of $\Fg$ in the linear space
$H$ iff there exists a set $\Fg_1,\ldots,\Fg_n$ of the Lie subalgebras of
$\Fg$, which form a dense connected composite and $T$ is its representation.
\enddefinition

Reducibility and irreducibility of representations of the Lie composites are
defined in the same manner as for Lie algebras. The set of representations of 
the fixed Lie composite is closed under the tensor product and, therefore, may 
be supplied by the structure of tensor category.

Let us consider two subalgebras $\Fp_{\pm}$ of the Witt algebra $\Fw^{\Bbb C}$ 
generated by $e_i$ with $i\ge-1$ and $i\le1$; note that $\Fp_+\cap\Fp_-=\sltwo$. 
The triple $(\Fw^{\Bbb C};\Fp_+,\Fp_-)$ is a dense connected Lie composite.

\proclaim{Proposition 3 [30]} The $q_R$--conformal symmetries in the Verma 
module $V_h$ over the Lie algebra $\sltwo$ form a representation of the Lie 
composite $(\Fw^{\Bbb C};\Fp_+,\Fp_-)$ and, hence, the composed representation
of the Witt algebra $\Fw^{\Bbb C}$.
\endproclaim

Theorems 1,2 and Proposition 3 imply Theorem 3.

\proclaim{Theorem 3} The operators $\nabla_X$ ($X\in\vir^{\Bbb C}$) in the 
Verma module $V_{h,c}$ over the Virasoro algebra form a representation
of the Lie composite $(\Fw^{\Bbb C};\Fp_+,\Fp_-)$ and, hence, the composed
representation of the Witt algebra $\Fw^{\Bbb C}$.
\endproclaim

\remark{Remark 3}
$$[L_X,\nabla_Y]=\nabla_{[X,Y]}$$
for any $X,Y\in\vir^{\Bbb C}$.
\endremark

It is naturally to consider the operators $\nabla_X$ as hidden symmetries in
the Verma modules $V_{h,c}$ over the Virasoro algebra and to try to set them 
free [31]. To do it one needs in some generalization of the procedure of the 
setting hidden symmetries free (as well as of definition of hidden symmetries)
in view of the infinite dimensionality of the family $\nabla_X$.

\definition{Definition 2}

{\bf A.} Let $\Fg$ be a $\Bbb Z$--graded Lie algebra and $\CA$ be a $\Bbb Z$--
graded associative algebra such that $\Fg$ is a $\Bbb Z$--graded subalgebra
of $\Der(\CA)$; a linear subspace $V$ of $\CA$ is called {\it a space of 
hidden symmetries\/} iff (1) $V$ is a $\Fg$--submodule of $\CA$, (2) the Weyl 
symmetrization defines mappings $W_{\pm}:S^{\cdot}(V_{\pm})\mapsto\CA$, 
and $W_0:S^{\cdot}(V_0)\mapsto\CA$, where $V_{\pm}$ and $V_0$ are the 
subspaces of $V$ of elements of positive, negative and zero degree, 
respectively (here $S^{\cdot}(H)$ is the symmetric algebra over a linear space 
$H$) such that the products $W_+(A)W_0(B)W_-(C)$ ($A\in S^{\cdot}(V_+)$,
$B\in S^{\cdot}(V_0)$, $C\in S^{\cdot}(V_-)$) are dense in $\CA$. The elements 
of $V$ are called {\it hidden symmetries with respect to $\Fg$}). A $\Bbb 
Z$--graded associative algebra $\CF$ such that $\Fg\subset\Der(\CF)$ is called 
{\it an algebra of the set free hidden symmetries\/} iff (1) $\CF$ is 
generated by $V$, (2) there exists a $\Fg$--equivariant epimorphism of 
algebras $\CF\mapsto\CA$, (3) the Weyl symmetrization define mappings 
$W_{\pm}:S^{\cdot}(V_{\pm})\mapsto\CF$ and $W_0:S^{\cdot}(V_0)\mapsto\CF$
so that their product $W_+W_0W_-$ defines the isomorphism of $\Fg$--modules
$S^{\cdot}(V_+)\otimes S^{\cdot}(V_0)\otimes S^{\cdot}(V_-)$ of $\Fg$--modules.

Let $\Fg$ be a Lie algebra, $V$ be a certain $\Fg$--module, $\CA_s$
be a family of associative algebras, parametrized by $s\in\CS$ such that
$\Fg\subset\Der(\CA_s)$, $\pi_s:V\mapsto\CA_s$ be a family of
$\Fg$--equivariant imbeddings such that $\pi_s(V)$ is a space of hidden
symmetries in $\CA_s$ with respect to $\Fg$ for a generic $s$ from $\CS$.  An
associative algebra $\CF$ is called {\it an algebra of the
$\CA_{s,s\in\CS}$--universally set free hidden symmetries\/} iff $\CF$ is an
algebra of the set free hidden symmetries corresponding to $V\simeq\pi_s(V)$
for generic $\CA_s$ ($s\in\CS$).
\enddefinition

The definitions are a natural $\Bbb Z$--graded generalizations of ones
from [31].

The following theorem, which unravel the algebraic structure of operators
$\nabla_X$, holds.

\proclaim{Theorem 4} The operators $\nabla_X$ in the Verma modules $V_{h,c}$
over the Virasoro algebra $\vir^{\Bbb C}$ form a family of hidden symmetries.
For any $c$ there exists a unique algebra $\CF_c$, which set these hidden
symmetries free $\End(V_{h,c})$--universally ($h\in\Bbb R$). The algebras
$\CF_c$ are the quotients of the central extension $\CF$ of the algebra
$\CF_0$ by the ideal generated by $\bold c-c\cdot\boldkey 1$, where
$\bold c$ is the central element.
\endproclaim

\remark{Remark 4} The theorem means that the algebra $\CF$ may be represented
in the model $M_c$ of the Verma modules $V_{h,c}$ over the Virasoro algebra
[27] for any $c$.
\endremark 

The theorem 4 easily follows from the dimension counting.

\head 3. Conclusions\endhead

Thus, some natural hidden symmetries in the Verma modules over the Virasoro
algebra are constructed in terms of geometric quantization. Their differential
geometric meaning is established and their expression via $q_R$--conformal
symmetries in the Verma modules over the Lie algebra $\sltwo$ is found. 
The analysis and the unraveling of the algebraic structure of these families 
of hidden symmetries are performed.

\Refs
\roster
\item"[1]" Juriev D., Infinite dimensional geometry of 
$M_1\!=\!\Diff_+(\Bbb S^1)/\PSL(2,\Bbb R)$ and $q_R$--conformal symmetries. 
E-print: math.RT/9806140.
\item"[2]" Juriev D.V., The vocabulary of geometry and harmonic analysis on
the infinite-di\-men\-sional manifold $\Diff_+(\Bbb S^1)/\Bbb S^1$. Adv.~Soviet 
Math. 2 (1991) 233-247.
\item"[3]" Kirillov A.A., Juriev D.V., Representations of the Virasoro
algebra by the orbit method. J.~Geom.~Phys. 5 (1988) 351-364 [reprinted
in: ``Geometry and physics. Essaya in honour of I.M.Gelfand. Eds. S.G.Gindikin
and S.G.Singer. Pitagora Editrice, Bologna and Elsevier Sci.Publ., Amsterdam,
1991].
\item"[4]" Zhelobenko D.P., Representations of the reductive Lie algebras
[in Russian]. Moscow, Nauka, 1993.
\item"[5]" Dixmier J., Alg\`ebres enveloppantes. Gauthier-Villars, Paris,
1974.
\item"[6]" Barut A., Raczka R., Theory of group representations and
applications. PWN -- Polish Sci.~Publ. Warszawa, 1977.
\item"[7]" Juriev D.V., Complex projective geometry and quantum projective
field theory. Theor. Math.~Phys. 101(3) (1994) 1387-1403.
\item"[8]" Berezin F.A., Quantization in complex symmetric spaces [in
Russian]. Izvestiya AN SSSR. Ser.~matem. 39(2) (1975) 363-402.
\item"[9]" Biedenharn L., Louck J., Angular momentum in quantum mechanics.
Theory and appli\-cations. Encycl.~Math.~Appl. V.8. Addison Wesley Publ.~Comp.
1981.
\item"[10]" Fuchs D.B., Cohomology of infinite-dimensional Lie algebras [in
Russian]. Moscow, Nau\-ka, 1984.
\item"[11]" Bott R., The characteristic classes of groups of diffeomorphisms.
Enseign.~Math. 23 (1977) 209-220.
\item"[12]" Kirillov A.A., Infinite-dimensional Lie groups, their invariants
and representations. Lect. Notes Math. 970 (1982) 101-123.
\item"[13]" Segal G., Unitary representations of some infinite-dimensional 
groups. Commun.~Math. Phys. 80 (1981) 301-342.
\item"[14]" Kirillov A.A., A K\"ahler structure on K-orbits of the group of
diffeomorphisms of the circle. Funct.Anal.Appl. 21 (1987) 322-324.
\item"[15]" Kirillov A.A., Juriev D.V., The K\"ahler geometry of the infinite
dimensional homogeneous space $M=\Diff_+(\Bbb S^1)/\Rot(\Bbb S^1)$.
Funct.~Anal.~Appl. 21 (1987) 284-293.
\item"[16]" Kirillov A.A., Elements of representation theory. Springer, 1976.
\item"[17]" Goluzin G.M., Geometric theory of functions of complex variables.
Amer.~Math.~Soc. 1968.
\item"[18]" Duren P.L., Univalent functions. Springer, 1983.
\item"[19]" Lehto O., Univalent functions and Teichm\"uller spaces. Springer,
1986.
\item"[20]" Kac V.G., Contravariant form for infinite-dimensional Lie algebras
and superalgebras. Lect.~Notes Phys. 94 (1979) 441-445.
\item"[21]" Feigin B.L., Fuchs D.B., Skew-symmetric invariant differential
operators on a line and Verma modules over Virasoro algebra [in Russian]. 
Funkt.~anal.~i ego prilozh. 16(2) (1982) 47-63.
\item"[22]" Feigin B.L., Fuchs D.B., Verma modules over Virasoro algebra [in
Russian]. Funkt.~anal.~i ego prilozh. 17(3) (1983) 91-92.
\item"[23]" Neretin Yu.A., Representations of Virasoro algebra and affine
algebras [in Russian]. Current Probl.~Math. Basic Directions. V.22, Moscow,
VINITI, 1988, P.163-224.
\item"[24]" Feigin B.L., Fuchs D.B., Representations of the Virasoro algebra.
In: ``Representations of inifinite dimensional Lie algebras''. Gordon and
Breach, 1991.
\item"[25]" Kostant B., Quantization and unitary representations. Lect.~Notes
Math. 170 (1970) 87-208.
\item"[26]" Kirillov A.A., Geometric quantization [in Russian]. Current
Probl.~Math. Basic Directions. V.4, Moscow, VINITI, 1985, P.141-178.
\item"[27]" Juriev D., A model of Verma modules over the Virasoro algebra,
St.-Petersburg Math.~J., 2 (1991) 401-417. 
\item"[28]" Adler S.L., Dashen R.F., Current algebras and applications to
particle physics. New York, 1968.
\item"[29]" Kobayashi Sh., Nomizu K., Foundations of differential geometry.
Intersci.~Publ., 1963/69.
\item"[30]" Juriev D., Topics in hidden symmetries. V. E-print: 
funct-an/9611003.
\item"[31]" Juriev D.V., An excursus into the inverse problem of
representation theory [in Russian]. Report RCMPI/95-04 (1995) [e-version:
mp\_arc:96-477].
\endroster
\endRefs
\enddocument